\documentclass[11pt]{elsarticle}

\usepackage{lineno}
\usepackage{amsmath, amssymb}
\usepackage{xspace, graphicx}
\modulolinenumbers[5]

\newcommand{\TT}{$T$-tetromino\xspace}
\newcommand{\TTs}{$T$-tetrominos\xspace}
\newtheorem{thm}{Theorem}
\newtheorem{lemma}[thm]{Lemma}

\begin{document}

\baselineskip=20pt
\begin{frontmatter}

\title{The Gap Number of the $T$-Tetromino}

\author{Robert Hochberg}\fnref{fn1}
\fntext[fn1]{972-721-5307 (w), 252-258-5284 (c)}
\address{Department of Mathematics, University of Dallas, Irving, TX}

\ead{hochberg@udallas.edu}

\begin{abstract}
A famous result of D. Walkup states that the only rectangles that may be tiled by the $T$-Tetromino are those in which both sides are a multiple of four. In this paper we examine the rest of the rectangles, asking how many $T$-tetrominos may be placed into those rectangles without overlap, or, equivalently, what is the least number of gaps that need to be present. 
We introduce a new technique for exploring such tilings, enabling us to answer this question for all rectangles, up to a small additive constant.
We also show that there is some number $G$ such that if both sides of the rectangle are at least 12, then no more than $G$ gaps will be required. We prove that $G$ is either 5, 6, 7 or 9.
\end{abstract}


\end{frontmatter}

\linenumbers

\section{Introduction and Definitions}
A well-known result of D. Walkup \cite{Walkup} is that the T-tetromino can tile a rectangle if and only if both sides of that rectangle are multiples of 4. This result has been extended in several ways. Mike Reid \cite{ReidKlarner} generalized the $T$-tetromimino to a class of $(8n-4)$-ominos, each of which tiles a rectangle if and only if both sides are multiples of 4, and other authors \cite{KornPak, Merino} have enumerated tilings of $4m\times 4n$ rectangles with $T$-tetrominos. 

Another direction of inquiry asks about those rectangles for which both sides are not multiples of 4. For odd $n=2k+1$, it seems possible that $k^2+k$ \TTs could be packed into the $n\times n$ square so that a single $1\times 1$ square remains uncovered. Wayne Goddard \cite{Goddard} showed
that this was not possible. Zhan \cite{Zhan} generalized this result to all $m\times n$ rectangles with $mn\equiv 1 \pmod 4$, showing again that it is impossible to leave just a single uncovered $1\times 1$ square. 

An uncovered $1\times 1$ square in a tiling may be thought of as a special $1\times 1$ tile, a {\it monomino}, employed as a second type of allowed tile. We will thus be interested in the following question: Given a rectangle, what is the greatest number of \TTs that may be placed therein when tiled with \TTs and monominos? Or, equivalently, what is the least number of monominos needed in such a tiling? Let us define 
$M(m, n)$ to be this number, the {\it gap number}. That is, $M(m, n)$ = the least number of monominos needed in a tiling of the $m\times n$ rectangle by \TTs and monominos. We observe that $M(m, n)\equiv mn \pmod 4$. 

We will show that $M(m, n)$ is unbounded as $m$ and $n$ range over the positive integers. It remains unbounded even if we fix $m$ to be 11. On the other hand, if $m$ and $n$ are both at least 12, then $M(m, n)\leq 9$. Thus, $G=\limsup_{m,n\rightarrow\infty}M(m, n)$ exists, and is at most 9. We call $G$ the {\it gap cap} of the \TT.

\paragraph*{Outline of the paper}In Section \ref{SecSmallOdd} we investigate rectangles of odd widths $w\leq 11$, for which $M(w, n)$ increases without bound as $n$ increases. Here we introduce the {\it fringe digraph} and show how the Bellman-Ford algorithm can be used to find lower bounds on the number of gaps needed, exactly matching constructive upper bounds. Section \ref{SecSmallEven} looks at even widths $w\leq 12$ and shows that for these widths, $M(w, n)$ remains bounded. Section \ref{Sec13and15} looks at rectangles of widths $w=13$ and $w=15$, and Section \ref{SecAllRectangles} shows that the gap cap of the \TT is at most 9.

\section{Rectangles of odd width $\leq 11$}\label{SecSmallOdd}

We consider rectangles of widths 1, 3, 5, 7, 9 and 11 in this section. We will cover width three in some detail, and then the others more quickly.

\subsection{Width 3}

\begin{thm}\label{Madeline}
\[
 M(3, n) =
  \begin{cases}
   5 & \text{if } n = 3 \\
   \lfloor n/3\rfloor       & \text{if } n\equiv 0\pmod 3, n\neq 3\\
   \lfloor n/3\rfloor + 3     & \text{if } n\equiv 1\pmod 3\\
   \lfloor n/3\rfloor + 2     & \text{if } n\equiv 2\pmod 3
  \end{cases}
\]
\end{thm}

The following lemma does most of the work for us and will have an analogue for width 9.

\begin{lemma}
In any tiling of a $3\times n$ rectangle by \TTs and monominos, each $3\times 3$ box must contain at least one monomino.
\end{lemma}

\begin{figure}[h]
\centering
\includegraphics[scale=0.9]{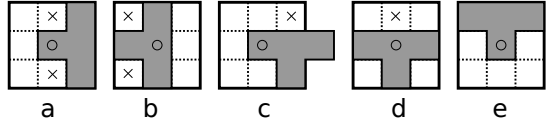}
\caption{Every $3\times 3$ box in any tiling of a $3\times n$ rectangle must have a gap.}
\label{FigWidth3}
\end{figure}

\paragraph{Proof of lemma:} 
Given a $3\times 3$ box consider how the center square is covered. If it is covered by a monomino, then we're done. If not, then it is covered by a \TT in some way, as shown (up to symmetry) in Figure \ref{FigWidth3}, where the dot represents the center square of some $3\times 3$ box. 

In cases a and b, it is impossible to cover both of the squares marked with ``X'' without using a monomino which will lie within that $3\times 3$ box. In cases c and d, the square marked with an ``X'' can not be covered with a \TT. In case e, the square below the dot may be covered with a \TT in only one way, up to symmetry, leaving two squares in that box that cannot both be tiled by \TTs. \hfill$\blacksquare$

\begin{figure}[h]
\centering
\includegraphics[scale=0.9]{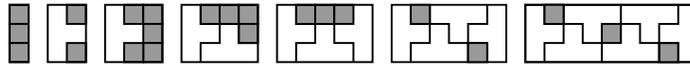}
\caption{Optimal $3\times n$ tilings for $n=1, 2, 3, 4, 5, 6, 9$.}
\label{FigWidth3SmallCases}
\end{figure}

\paragraph{Proof of Theorem:} The cases for $n\leq 5$ are shown in Figure \ref{FigWidth3SmallCases}. In what follows, we assume $n\geq 6$. Let $n=3k + r$ where $0\leq r\leq 2$. Then within a $3\times n$ rectangle we can find $k$ disjoint $3\times 3$ blocks, each of which contains a monomino, so that 

\begin{equation}\label{EQ3nlowerbound}
M(3, n) \geq k =\lfloor n/3\rfloor.
\end{equation}

For $r=0$, the $3\times 6$ and $3\times 9$
rectangles in Figure \ref{FigWidth3SmallCases} may be combined to show that $M(3, n) \leq\lfloor n/3\rfloor$, finishing the proof for the case where $n\equiv 0\pmod 3$. 

Now suppose $r=1$. Since the $3\times (3k+1)$ rectangle contains three more squares than the $3\times 3k$ rectangle, and recalling that $M(m, n)\equiv mn\pmod 4$ we have that the number of gaps required for the $3\times (3k+1)$ rectangle must equal $k+3+4t$ for some integer $t$. By (\ref{EQ3nlowerbound}), $t$ can't be negative, giving $M(3, n)\geq k+3$. This number of gaps can be realized by simply appending a column of three monominos to an optimally-tiled $3\times 3k$ rectangle, finishing the proof when $n\equiv 1\pmod 3$. If $r=2$, then the number of gaps must be $k+6+4t$ for some integer $t$. Again, by (\ref{EQ3nlowerbound}), $t\geq -1$, giving $M(3, n)\geq k+2$. This bound may be achieved by appending to the optimally-tiled $3\times 3k$ rectangle two columns containing one \TT and two monominos.\hfill$\blacksquare$

\subsection{Rectangles of Width 9}
We treat width 9 next because it is resolved in much the same way as width 3. We have a similar lemma:

\begin{lemma}
In any tiling of a $9\times n$ rectangle by \TTs and monominos, each $9\times 17$ box must contain at least one monomino.
\end{lemma}

The proof of this lemma was carried out by exhaustive computer search, where we discovered that it is possible to completely cover by \TTs alone sixteen consecutive columns in a width-9 strip, but not seventeen. Thus $M(9,n)\geq \lfloor n/17 \rfloor$.

\begin{figure}[h]
\centering
\includegraphics[angle=0, scale=2]{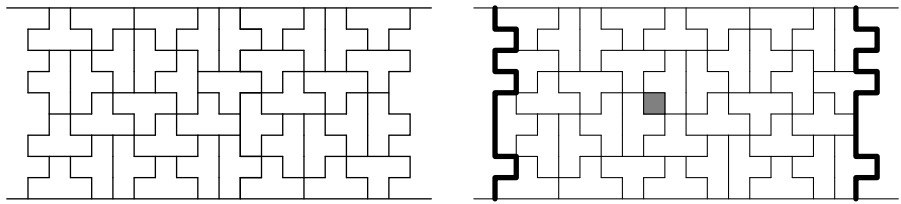}
\caption{Left: Sixteen columns of a width-9 strip tiled without gaps. Right: A $9\times 17$ cylinder with one monomino.}
\label{FigCylinder9}
\end{figure}

On the other hand, we have the tiling fragment shown in Figure \ref{FigCylinder9} with which we may periodically cover arbitrarily many columns with only one monomino every seventeen columns. We have therefore proven:

\begin{thm}
$M(9, n) = \lfloor n/17\rfloor + e_9(n)$ for some bounded function $e_9(n)$.\hfill$\blacksquare$
\end{thm}

We conjecture that $e_9(n) \le 6$ for $n\geq 2$.

\subsection{Rectangles of Widths 1, 5, 7 and 11, and the Fringe Digraph}

We ask the reader to verify that $M(1, n) = n$.

\begin{figure}[h]
\centering
\includegraphics[scale=1.4]{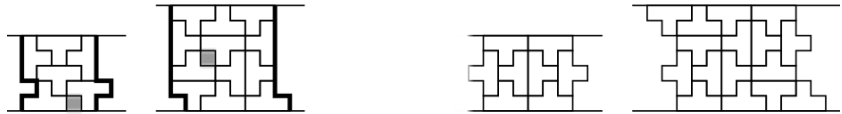} 
\caption{Widths 5 and 7. Left: $5\times 5$ and $7\times 7$ cylinders with one monomino. Right: It is possible to tile many columns without monominos.}
\label{FigWidths5and7}
\end{figure}

The left side of Figure \ref{FigWidths5and7} shows a fragment for tiling width-5 rectangles periodically, using one monomino for each five columns, giving $M(5, n) \leq n/5 + C$ for some constant C. The second fragment gives $M(7, n) \leq n/7 + C$. Unfortunately, the method used for widths 3 and 9 will not work to establish corresponding lower bounds for $M(5, n)$ and $M(7, n)$. When tiling a strip of width 5, it is possible to cover six consecutive columns without a monomino, but not seven, which would give the lower bound $M(5, n) \geq n/7$. And for a strip of width 7, we may tile up to eight columns without a monomino, but not nine, giving the lower bound $M(7, n)\geq n/9$. (See the right side of Figure \ref{FigWidths5and7}.) This leaves a gap between our upper and lower bounds.

\begin{figure}[h]
\centering
\includegraphics[scale=0.3]{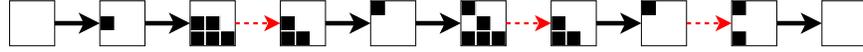} 
\caption{The sequence of fringes in our $3\times 9$ tiling from Figure \ref{FigWidth3SmallCases}.}
\label{FigFringe3x9}
\end{figure}

If a greedy tiler were to try to tile a rectangle of width 5 or 7 using one of the fragments given on the right in Figure \ref{FigWidths5and7}, he would discover that many monominos would be needed nearby, ``making up'' for the monomino-free columns. The case-by-case analysis required to quantify this rigorously quickly becomes onerous. We therefore turn again to the computer, and make use of the {\it fringe digraph}, which we now describe. Suppose we tile a rectangle with \TTs and monominos, starting from the left side of the rectangle, and placing tiles one-by-one so that we always cover the lowest of the left-most untiled squares. Then since our tiles both fit within a $2\times 3$ box, our tiling can never have more then three partially-tiled columns, called the {\it fringe}. For example, the tiling of the $3\times 9$ rectangle shown in Figure \ref{FigWidth3SmallCases} yields the sequence of fringes shown in Figure \ref{FigFringe3x9}.

\begin{figure}[h]
\centering
\includegraphics[scale=0.2]{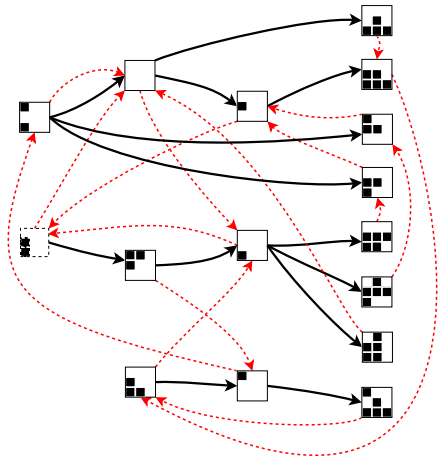}
\caption{The fringe digraph $F_3$.}
\label{FigFringeDigraph3}
\end{figure}

When we begin tiling a rectangle the fringe is empty, and this gives us the first node in the fringe digraph. The rest of the digraph is constructed recursively: Each time we discover a new fringe we create a corresponding node $N$ in the digraph, and consider all ways to place a tile that covers the lowest of its left-most open squares. If it yields a previously-seen fringe, we find the node $C$ corresponding to that fringe, and add a directed edge from $N$ to $C$. If it gives a new fringe, we recursively visit that fringe, and then add the arc from $N$ to that new node. Let $F_w$ denote the fringe digraph for width $w$. $F_3$ is shown in Figure \ref{FigFringeDigraph3}.

We distinguish between those edges that correspond to adding a \TT and those that correspond to adding a monomino. Let us call them $T$-edges and $M$-edges respectively. In Figure \ref{FigFringeDigraph3}, the $M$-edges are thinner and dotted. Not surprisingly, this digraph contains much information. For example, we can get a new proof that $M(3, n)\leq n/3 + C$ by considering the 3-cycle on the three lowest nodes of that figure: That cycle contains two $T$-edges and one $M$-edge, accounting for 4+4+1=9 squares, or three columns, containing a single monomino.

The fringe digraph $F_w$ can also be used to find lower bounds for $M(w, n)$. Suppose, for example, we put weights on the edges of $F_3$: $-1$ on $T$-edges and 2 on $M$-edges. It is not hard to verify that every cycle in the resulting weighted digraph has non-negative weight. In particular, any cycle starting and ending at the node corresponding to an empty fringe will have non-negative weight, and these are the cycles that correspond to tilings of rectangles. But this weight can be non-negative only if there is at least one $M$-edge for every two $T$-edges, that is, if there is at least one monomino for every three columns. 

We ran this same analysis for the fringe digraphs $F_5, F_7, F_9$ and $F_{11}$ using the Bellman-Ford algorithm to discover negative-weight cycles. The results are shown in Table \ref{TableFringe}. Let $-1$ and $m$ be the weights assigned to each $T$-edge and $M$-edge respectively.
 Then the absence of negative-weight cycles implies that any cycle containing $\tau$ $T$-edges and $\mu$ $M$-edges must satisfy $-\tau + m\mu \geq 0$, or equivalently, $\tau \leq\mu m$. Such a cycle would correspond to covering $4\tau + \mu$ squares of the rectangle, or $(4\tau + \mu)/w$ columns if the rectangle had width $w$. That is, the ratio of columns to monominos is $(4\tau + \mu)/\mu w$. Combined with the previous inequality, we have

\begin{equation}
\frac{\mbox{\# columns}}{\mbox{\# monominos}}=\frac{4\tau + \mu}{\mu w}
\leq \frac{4m\mu + \mu}{\mu w} = \frac{4m+1}{w}.
\end{equation}

\begin{figure}[h]
\centering
\begin{tabular}{|c|c|c|c|r|}
\hline
Fringe Digraph $F_w$&$m$		&$t$ 	& $\frac{4m+1}{w}$ & \# nodes in $F_w$\\
\hline
\hline
$F_3$	&2		&$-1$	&3	&16\\
\hline
$F_5$	&6		&$-1$	&5	&182\\
\hline
$F_7$	&12		&$-1$	&7	&1757\\
\hline
$F_9$	&38		&$-1$	&17	&15,496\\
\hline
$F_{11}$	&96		&$-1$	&35	&129,500\\
\hline
\end{tabular}
\caption{If we assign weight $m$ to the $M$-edges and weight $t$ to the $T$-edges of the fringe digraphs, then there will be no negative-weight cycles, so that there must be at least one monomino for every $(4m+1)/w$ columns.}
\label{TableFringe}
\end{figure}

By selecting appropriate (that is, the smallest) values for $m$ so that $F_w$ has no negative-weight cycles, as shown in the table, we obtain our desired lower bounds: In any tiling of a width-5 rectangle, there must be at least one monomino for every five columns, and for width-7, at least one monomino for every seven columns.

\begin{figure}[h]
\centering
\includegraphics[scale=1.3]{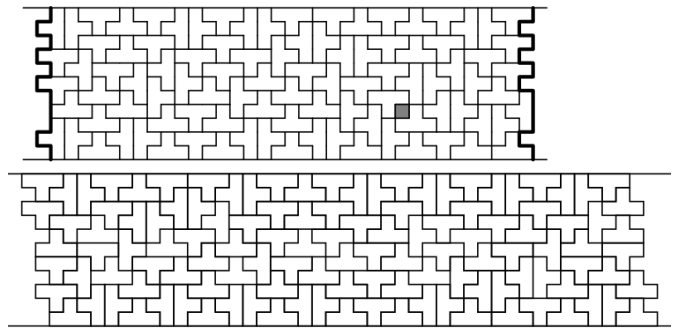} 
\caption{Top: A tiling fragment of width 11 that may tile periodically with one monomino in every 35 columns. Bottom: 42 columns tiled without any monominos.}
\label{FigWidth11}
\end{figure}

Finally, we consider rectangles of width 11. As the table proves, in any tiling of an $11\times n$ rectangle we will find at least $n/35$ monominos. The tiling fragment at the top of Figure \ref{FigWidth11} shows that this bound is obtainable. This bound surprised us a bit, after having first found that 42 columns could be tiled without any gaps, as shown at the bottom of that figure.

We therefore have the following theorem:

\begin{thm}
$M(5, n) = \lfloor n/5\rfloor + e_5(n)$, 
$M(7, n) = \lfloor n/7\rfloor + e_7(n)$,
$M(11, n) = \lfloor n/35\rfloor + e_{11}(n)$ for some bounded functions $e_5(n)$, $e_7(n)$ and $e_{11}(n)$.\hfill$\blacksquare$
\end{thm}

We conjecture that $e_5(n)\leq 4$ for $n>15$, $e_7(n)\leq 4$ for $n>10$ and $e_{11}(n)\leq 7$ for $n > 1$.

\section{Rectangles of Even Width at most 12}\label{SecSmallEven}

The situation is very different for even-width rectangles. Here we will see that for even $w$, $M(w, n)$ remains bounded as $n$ grows. Walkup's theorem tells us that $M(w, n)$ cannot be zero unless both $w$ and $n$ are multiples of 4. Together with the observation that $M(w, n)\equiv wn\pmod 4$, this gives the value of $M(w, n)$ for most rectangles of even width $w$. 

\begin{thm}\label{ThmEvenw}
For even $w\leq 12$ and any $n\geq 2$, the value of $M(w, n)$ is the least positive integer $k$ allowed by Walkup's theorem, such that $k\equiv wn\pmod 4$, with only the following exceptions: $M(10, n) = 6$ for $n = 3, 7, 11, 15, 19$ and $23$.
\end{thm}

\subsection{Proof of Theorem \ref{ThmEvenw}}
\begin{figure}[h]
\centering
\includegraphics[scale=1.7]{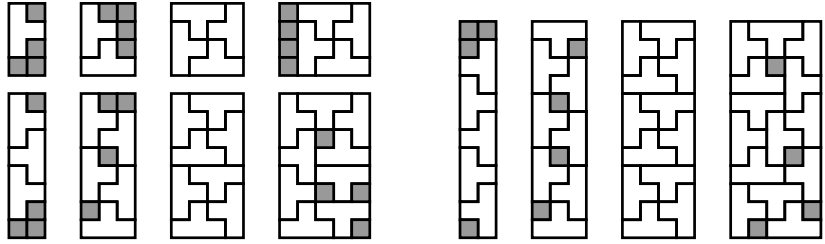} 
\caption{Rectangles of width 4, 8 and 12 require 0 or 4 monominos.}
\label{FigWidth4-8-12}
\end{figure}

The proof for width 2 is left to the reader. Rectangles of width 4, 8 and 12 require 0 or 4 monominos for length at least 2. (We note that $M(n, 1) = n$.) Figure \ref{FigWidth4-8-12} shows optimally-tiled rectangles of lengths 2, 3, 4 and 5, and to these may be concatenated any number of the $w\times 4$, monomino-free rectangles, to obtain the theorem for widths 4, 8 and 12.

\begin{figure}[h]
\centering
\includegraphics[scale=1.7]{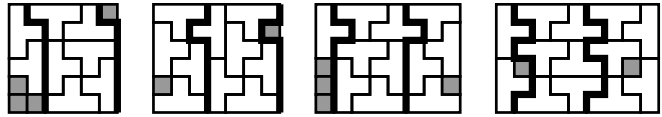}
\caption{Rectangles of width 6 require 2 or 4 monominos. The central section of each rectangle may be deleted or replicated, to obtain $6\times n$ rectangles for all $n\geq 2$.}
\label{FigWidth6}
\end{figure}

Rectangles of width 6 require two monominos when $n$ is odd, and four when $n$ is even. The four rectangles shown in Figure \ref{FigWidth6} have lengths in each of the congruence classes (mod 4) and central regions of width 4 that may be deleted or repeated to obtain tilings of $6\times n$ rectangles for all $n\geq 2$, matching the bounds of the theorem.

\begin{figure}[h]
\centering
\includegraphics[scale=1.05]{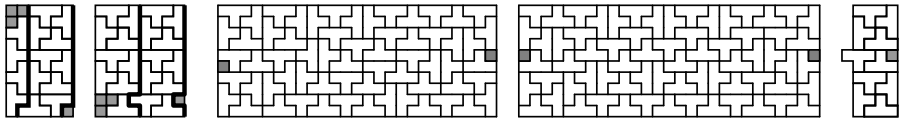} 
\caption{The central section of the left two figures may be deleted or repeated, showing that $M(10, k) = 4$ for all even $k\geq 2$. The two large rectangles show $M(10, 25) = M(10, 27) = 2$, and in each of those two rectangles, the rightmost monomino may be replaced by the fragment on the right to show that $M(10, k) = 2$ for all odd $k\geq 25$.}
\label{FigWidth10}
\end{figure}

Rectangles of width 10 are handled in a similar fashion to those of width 6. The two rectangles on the left of Figure \ref{FigWidth10} show that $M(10, k) = 4$ for all even $k\geq 2$, again making use of monomino-free central sections of width 4 that may be deleted or repeated. For odd $k\geq 25$ Figure \ref{FigWidth10} shows $10\times 25$ and $10\times 27$ rectangles, each of which may be extended by any number of copies of the four-column unit shown. The odd values of $k$ from 3 to 23 can then be verified by computer, which we have done. This completes the proof of Theorem \ref{ThmEvenw}.\hfill$\blacksquare$

\section{Rectangles of Width 13 and Width 15}\label{Sec13and15}

\begin{figure}[h]
\centering
\includegraphics[scale=1.6]{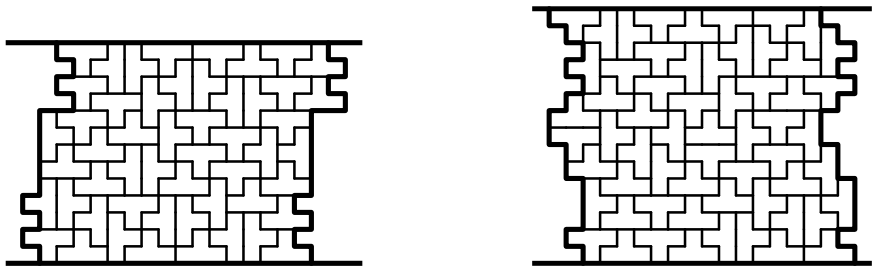} 
\caption{Cylinders of widths 13 and 15.}
\label{FigCyl13and15}
\end{figure}

After having found the rate of growth of $M(w, n)$ for the small, odd $m$ discussed above, it came as a great surprise to find that for width 13, it was possible to tile arbitrarily many columns without needing a monomino. The tiling fragment on the left in Figure \ref{FigCyl13and15} shows how a periodic tiling of a strip of width 13 may be obtained with \TTs only. We call this fragment a {\it cylinder} of size $13\times 16$, and it is the smallest odd-width cylinder tileable by \TTs alone. A $15\times 16$ cylinder is shown on the right in that Figure.

Walkup's theorem says that $M(w, n)$ can not be 0 for these widths, and Zhan's theorem \cite{Zhan} says that it also can not be 1. 

\begin{thm}\label{ThmWidth13and15}
For $w=13$ and $w=15$ and any $n\geq 2$, the value of $M(w, n)$ is the least positive integer $k\geq 2$ such that $k\equiv wn\pmod 4$, with only the following exceptions: $M(13, 3) = 7$, $M(13, n) = 6$ for $n$=14, 18, 22, 30, 34, 38, $M(15, 5) = 7$ and $M(15, n)=6$ for $n$=10, 14, 18, 22, 26. 
\end{thm}

\begin{figure}[h]
\centering
\includegraphics[scale=0.5]{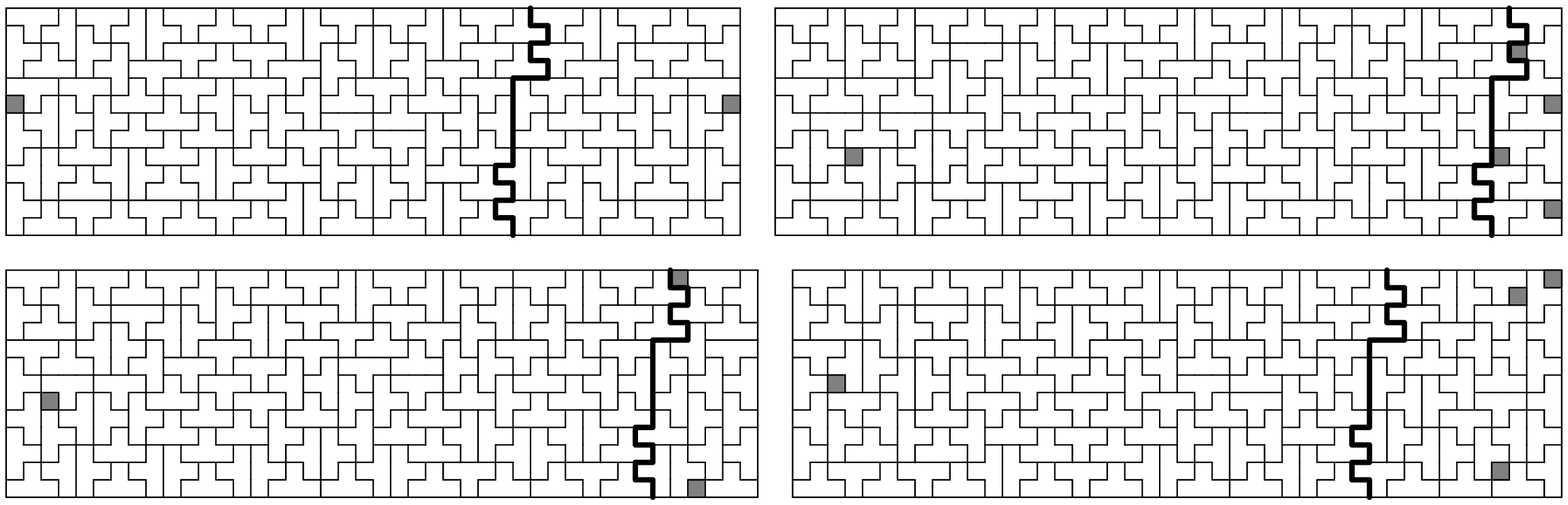}
\caption{Rectangles of width 13 which may be mixed and matched at the border, together with the cylinder in Figure \ref{FigCyl13and15}, to obtain all optimally-tiled $13\times n$ rectangles for $n\geq 39$.}
\label{FigWidth13}
\end{figure}

\paragraph{Proof of Theorem}
Figure \ref{FigWidth13} shows four rectangles of width 13. These rectangles all share a common boundary, as shown, allowing the left and right parts to be mixed and matched to create sixteen different rectangles, each of which is tiled with the number of monominos guaranteed by the theorem. The lengths of these sixteen rectangles lie in the sixteen different congruence classes (mod 16), as shown in the middle column of the table in Figure \ref{TableLengths}. Furthermore, this boundary is shared by the $13\times 16$ cylinder shown in Figure \ref{FigCyl13and15}, allowing us to extend those rectangles to obtain all $13\times n$ rectangles for $n\geq 39$, and some smaller ones. The few remaining, smaller rectangles may be checked by computer, which we have done.

\begin{figure}[h]
\centering
\begin{tabular}{|c|c|c|}
\hline
mod 16  &Width 13 start		&Width 15 start\\
\hline
\hline
0&48&16\\
\hline
1&33&33\\
\hline
2&50&50\\
\hline
3&35&19\\
\hline
4&52&20\\
\hline
5&37&37\\
\hline
6&54&38\\
\hline
7&39&23\\
\hline
8&40&24\\
\hline
9&41&25\\
\hline
10&42&42\\
\hline
11&43&27\\
\hline
12&44&28\\
\hline
13&45&29\\
\hline
14&46&46\\
\hline
15&47&15\\
\hline
\end{tabular}
\caption{Lengths that may be obtained for rectangles of width 13 and 15 by mixing and matching the right and left halves of the rectangle shown in Figures \ref{FigWidth13} and \ref{FigWidth15}.}
\label{TableLengths} 
\end{figure}

\begin{figure}[h]
\centering
\includegraphics[scale=0.5]{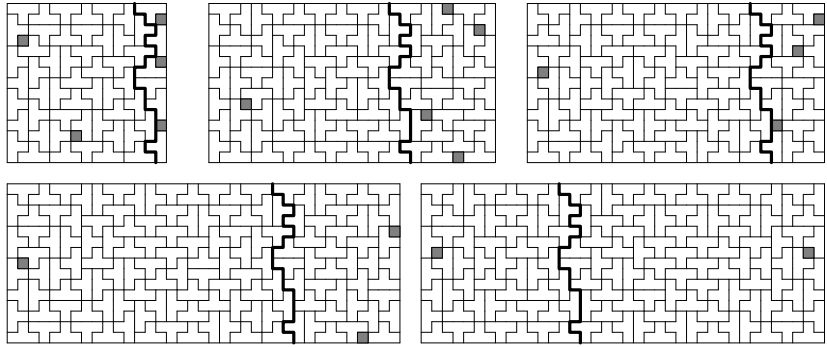} 
\caption{Rectangles of width 15 which may be mixed and matched at the border, together with the cylinder in Figure \ref{FigCyl13and15}, to obtain all optimally-tiled $15\times n$ rectangles for $n\geq 37$.}
\label{FigWidth15}
\end{figure}

In the same fashion, we may mix and match the right and left parts of the rectangles shown in Figure \ref{FigWidth15} to obtain rectangles with the lengths shown in Figure \ref{TableLengths} for width 15, hitting all sixteen congruence classes (mod 16) and tiled with the number of monominos guaranteed by the theorem. Again, the boundary shown is shared by the $15\times 16$ cylinder shown in Figure \ref{FigCyl13and15}, yielding optimally-tiled rectangles for all lengths at least 37. The few remaining, smaller rectangles may be checked by computer, and we have done that. This finishes the proof of Theorem \ref{ThmWidth13and15}.\hfill$\blacksquare$

\section{Rectangles with Both Sides at Least 12}\label{SecAllRectangles}
We may now find an upper bound on the gap cap ($=\limsup_{m,n\rightarrow\infty}M(m, n)$) of the \TT.

\begin{thm}\label{ThmGapNumber}
If $m$ and $n$ are both at least 12, then $M(m, n)\leq 9$. 
\end{thm}

\begin{figure}[h]
\centering
\includegraphics[scale=0.55]{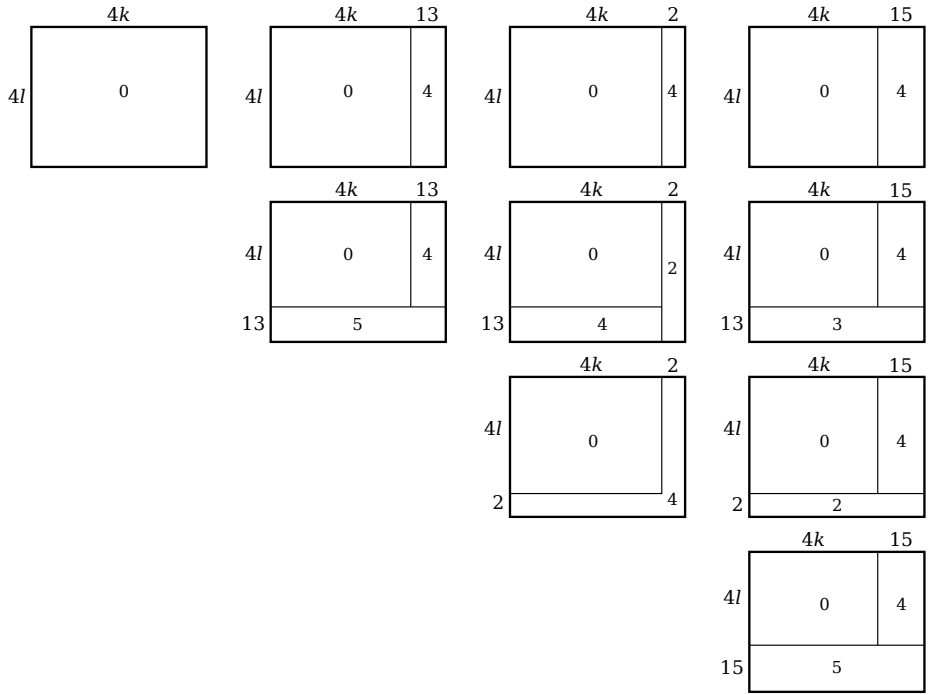} 
\caption{Decomposing rectangles into pieces that require few gaps.}
\label{FigStripChart}
\end{figure}

{\it Proof:} Theorem \ref{ThmEvenw} suffices to prove our theorem for rectangles that have at least one side of length exactly 12. We therefore assume both $m$ and $n$ are at least 13, and write $m=4l+r_1$ and $n=4k+r_2$ where $r_1$ and $r_2$ are in $\{0, 13, 2, 15\}$. We cut strips with width $r_1$ and/or $r_2$ off the bottom and/or right side of the rectangle, so that a $4l\times 4k$ block remains. This block requires no monominos, and the numbers of monominos needed by the strips of width 2, 13 or 15 are given by previous theorems. These numbers are shown inside the rectangles in Figure \ref{FigStripChart}. Note that the bottom rectangle in the third column, the case where $m, n\equiv 2\pmod{4}$, is different than the others. In that case an L-shaped border may be removed from the right and bottom, and it is easily shown that this may be tiled with just four monominos. In all cases, the total number of monominos needed is not more than nine, and the theorem is proved.\hfill$\blacksquare$

\subsection{Discussion}
We conjecture that the actual gap cap of the \TT is five rather than nine. There are some small rectangles of widths 13 and 15 that require six monominos, as listed in the statement of Theorem \ref{ThmWidth13and15}, so that the threshold would have to be increased from 12 to 39 if this conjecture were true. 

It would be nice to find the exact values of the functions $e_5(n)$, $e_7(n)$, $e_9(n)$ and $e_{11}(n)$ mentioned in Section \ref{SecSmallOdd}. Even better would be to find a way to ``turn the corner'' in the $(4k+13)\times (4l+13)$ and $(4k+15)\times (4l+15)$ cases, and show that by combining the right and bottom rectangle slices, only five monominos are needed. This would reduce our bound on the gap cap to seven. Perhaps doing something similar in the other cases as well could show that the gap cap is indeed five.

Finally, we note that the gap cap for some other small polyominos may be found without much difficulty, as shown in Figure \ref{FigSmallGaps}. We conjecture that the gap cap of the $S$-tetromino is also unbounded. And except for the few cylinders mentioned in this paper, we have not explored tilings of other surfaces, though these questions seem to be interesting as well.

\begin{figure}[h]
\centering
\begin{tabular}{|c|c|c|}
\hline
Polyomino  &Gap Cap		&Threshold\\
\hline
\hline
Monomino	& 0 & 1\\
\hline
Domino & 1 &1\\
\hline
$I$-tromino & 2 & 3\\
\hline
$L$-tromino & 2 & 4\\
\hline
$O$-tetromino& $\infty$&\\
\hline
$L$-tetromino & 4 & 4\\
\hline
$I$-tetromino & 4 & 4\\
\hline

\end{tabular}
\caption{Gap caps and thresholds for some small polyominos.}
\label{FigSmallGaps} 
\end{figure}

\section{Acknowledgments}
The author would like to thank Madeline Sargent who began this research and proved Theorem \ref{Madeline} while a high school student, and Zachary Dain who as an undergraduate first discovered that the \TT could tile 20,000 columns of a width-13 strip without gaps.

\bibliography{GapNumber}

\end{document}